\documentclass[11pt, a4paper]{article}

\usepackage[T1]{fontenc}
\usepackage{lmodern}
\usepackage[margin=25mm]{geometry}
\usepackage{amsmath, amssymb, amsthm}
\usepackage{microtype}

\usepackage{tikz}
\usetikzlibrary{decorations.pathreplacing,calc}
\usepackage{float}
\usepackage[font=small, labelfont=bf, skip=8pt]{caption}
\usepackage[hidelinks]{hyperref}
\hypersetup{
	pdftitle={Weak oddness to resistance ratios in snarks}
}

\AtBeginDocument{%
	\setlength{\abovedisplayskip}{8pt plus 2pt minus 2pt}
	\setlength{\belowdisplayskip}{8pt plus 2pt minus 2pt}
	\setlength{\abovedisplayshortskip}{4pt plus 2pt minus 1pt}
	\setlength{\belowdisplayshortskip}{6pt plus 2pt minus 1pt}
}

\newtheorem{theorem}{Theorem}

\newtheorem{conjecture}{Conjecture}
\newtheorem{lemma}{Lemma}
\newtheorem{problem}{Problem}

\definecolor{one}{RGB}{35, 84, 160}
\definecolor{two}{RGB}{175, 48, 43}
\definecolor{three}{RGB}{26, 114, 77}

\tikzset{
	vertex/.style={circle, fill=black, inner sep=1.7pt},
	c0/.style={black, very thick, densely dashed},
	c1/.style={one, thick},
	c2/.style={two, thick},
	c3/.style={three, very thick},
	el/.style={fill=white, inner sep=1pt, font=\small}
}

\allowdisplaybreaks

\title{Weak oddness to resistance ratios in snarks}

\author{	Imran Allie \\	
	\small Department of Mathematics and Applied Mathematics\\[-0.8ex]
	\small University of Cape Town\\[-0.8ex]
	\small Rondebosch, Cape Town, 7700, South Africa\\
	\small \texttt{imran.allie@uct.ac.za}}

\begin{document}
	
	\maketitle
	
	\begin{abstract}
		Let $G$ be a bridgeless cubic graph. If $G$ is not 3-edge-colourable, we call it a \textit{snark}. The \textit{resistance} of $G$, denoted $r(G)$, is the minimum number of edges which can be removed from $G$ in order to render $G$  3-edge-colourable. The \textit{oddness} of $G$, denoted $\omega(G)$, is the minimum number of odd components in any 2-factor of $G$. The \textit{weak oddness} of $G$, denoted $\omega'(G)$, is the minimum number of odd components in any even factor of $G$. In their comprehensive survey paper on measurements of uncolourability in cubic graphs, Fiol et al. conjectured the following two statements for a snark $G$: $\omega(G) \leq 2r(G)$ and $\omega'(G) \leq 2r(G)$ \cite{FiolMazzuoccoloSteffen2018}. Previously, we disproved the former by proving that the ratio of oddness to weak oddness, and hence the ratio of oddness to resistance, can in fact be arbitrarily large \cite{Allie2019}. In this paper, we prove a slightly stronger version of the latter: $\omega'(G) < 2r(G)$. The proof centres around substructures in $G$, given a minimal colouring $f$, which we call \textit{meta-Kempe chains}. 
	\end{abstract}
	
	
	\section{Introduction}
	
	\textit{Snarks} are bridgeless cubic graphs whose edges cannot be properly coloured with three colours. The significance of snarks derives mainly from the fact that they are the only class of graphs which may possibly contain counterexamples to many major open conjectures in graph theory. The cycle double cover conjecture (notwithstanding OpenAI's proof from mid 2026), Tutte's $5$-flow conjecture, and Fulkerson's conjecture are some examples of these major open conjectures. In particular, these potential counter-examples have also been shown, in some cases, to necessarily have girth greater than 4 and cyclic connectivity greater than 3. Such snarks are often called `non-trivial' in the literature. Those snarks which are otherwise `trivial' are easily reducible to smaller snarks by well-defined reductions. We do, however, believe that the idea of what constitutes triviality in snarks should be considered more deeply. See \cite{Allie2022} for further ideas on this discussion. In this paper we use the simple convention of snarks as uncolourable bridgeless cubic graphs, without concerning ourselves with connectivity or girth.
	
	Researchers have endeavoured to understand snarks by defining and investigating parameters which `measure', in various ways, how far a snark is from being 3-edge-colourable. Such parameters have been termed \textit{measures of uncolourability}. In their survey on measures of edge-uncolourability of cubic graphs,
	Fiol, Mazzuoccolo and Steffen consider the relationship between many such parameters, and put forth many problems and conjectures \cite{FiolMazzuoccoloSteffen2018}. Let $G$ be a bridgeless cubic graph. If $G$ is not 3-edge-colourable, we call it a \textit{snark}. The \textit{resistance} of $G$, denoted $r(G)$, is the minimum number of edges which can be removed from $G$ in order to render $G$ 3-edge-colourable. The \textit{flow resistance} of $G$, denoted $r_f(G)$, is the minimum number of zero edges in a $4$-flow of $G$. Fiol et al. conjectured that  $r_f(G) \leq r(G)$. This was disproved in \cite{AllieMacajovaSkoviera2022}. Further, in \cite{AllieMacajovaSkoviera2024}, it was shown that the ratio difference of flow resistance to resistance can be arbitrarily large. The \textit{colouring defect} of $G$, denoted $\mu_3(G)$, is the minimum number of edges not covered by $3$ perfect matchings of $G$. Recent results involving this measure relating to others include \cite{Allie2026, KarabasEtAl2022, KarabasEtAl2024, KarabasEtAl2025}. In \cite{Allie2026} we showed that snarks exist, essentially, for any given value of resistance, oddness and colouring defect. Again showing that an arbitrarily large ratio difference is possible between said parameters.
	
	The \textit{oddness} of $G$, denoted $\omega(G)$, is the minimum number of odd components in any 2-factor of $G$. The \textit{weak oddness} of $G$, denoted $\omega'(G)$, is the minimum number of odd components in any even factor of $G$. Oddness and weak oddness have sometimes been used interchangeably in the literature, with their equality implicitly assumed. This is alluded to in \cite{FiolMazzuoccoloSteffen2018}. Luko\v{t}ka and Maz\'ak showed that these parameters are distinct by constructing a bridgeless cubic graph with resistance $12$, weak oddness $14$, and oddness $16$ in \cite{LukotkaMazak2018}. They further showed that the arithmetic difference between oddness and weak oddness can be arbitrarily large. This implies as well that the arithmetic difference between oddness and resistance can be arbitrarily large. On the other hand, Fiol et al. proposed that the ratio difference is bounded. They conjectured the following two well-motivated statements for every bridgeless cubic graph $G$:
	
	\begin{conjecture} {\rm \cite{FiolMazzuoccoloSteffen2018} \label{con_oddness_resistance}
			Let $G$ be a bridgeless cubic graph. Then $\omega(G)\leq 2r(G)$.}	
	\end{conjecture}
	
	\begin{conjecture} {\rm \cite{FiolMazzuoccoloSteffen2018} \label{con_weakoddness_resistance}
			Let $G$ be a bridgeless cubic graph. Then $\omega'(G)\leq 2r(G)$.}	
	\end{conjecture}
	
	We disproved Conjecture \ref{con_oddness_resistance} in \cite{Allie2019}, where bridgeless cubic graphs of resistance three, weak oddness four, and arbitrarily large oddness are constructed. The construction forces the traversal of a cycle in any 2-factor through an arbitrarily large number of Petersen graph fragments, thereby forcing the existence of many odd cycles in the 2-factor, whilst keeping resistance constant. The allowance of isolated vertices in an even factor makes it impossible to force such a traversal. In \cite{FiolMazzuoccoloSteffen2018}, the authors prove that if $r(G)=3$, then $\omega'(G)=4$. The method involves Kempe chains. Given a proper edge-colouring, a \textit{Kempe chain} is a path such that edges are coloured alternately with two colours. The proof considers Kempe chains initiating at `conflicting edges' with respect to the colouring and how they possibly and/or necessarily intersect each other, then appropriately swapping colours on subchains to reveal an even factor with at most four odd components. We believe this method can be reproduced for snarks with $r(G)=4$ and $r(G)=5$, to show that $\omega'(G) \leq 6$ and $\omega'(G) \leq 8$, respectively, but proves insufficient for $r(G) \geq 6$. In this paper, we define \textit{meta-Kempe chains} in proving that $\omega'(G) < 2r(G)$ for any snark $G$.  
	
	\section{Preliminaries}
	
	In this paper, we only consider simple graphs (no parallel edges or loops) and we assume the reader to have knowledge of standard graph theoretical definitions. We do present our conventions relating to proper edge-colourings of cubic graphs. Throughout, let $G$ be a simple cubic graph. An \emph{edge-colouring}, $f$, assigns a colour to each edge. It is \emph{proper} if adjacent edges receive different colours. A \emph{colour class} is the set of edges receiving a given colour. A graph is \emph{$k$-edge-colourable} if it has a proper edge-colouring using at	most $k$ colours. The \emph{resistance} of $G$, $r(G)$, is the minimum number of edges whose deletion makes $G$ $3$-edge-colourable. For cubic graphs, this is also the minimum possible size of a colour class in a proper $4$-edge-colouring. Throughout, we use the colour set $\{0,1,2,3\}$ with $|f^{-1}(0)|=r(G)$. Such a colouring $f$ is called a \textit{minimal colouring}. Given $G$ and $f$, an \textit{$(a,b)$-Kempe chain} is a path in $G$ such that edges are coloured alternately with colours $a$ and $b$. Where context is clear, we will simply say \textit{Kempe chain} or \textit{$(a,b)$-chain}. 
	
	An \textit{$(a,b)(a,c)$-meta-Kempe chain} is a walk  which is an alternating sequence of $(a,b)$-chains and $(a,c)$-chains, with the terminal vertex of a given $(a,b)$-chain being the initial vertex in the next $(a,c)$-chain in the sequence (or vice versa). A vertex which is either a terminal vertex on some chain or an initial vertex on some chain is called a \textit{switching} vertex. We will always consider the terminal edge on a subchain to be the edge coloured $b$ or $c$, not $a$. Note that the proceeding edge immediately after the switching vertex could be coloured $a$ or not $a$. Each maximal Kempe chain between two switching vertices is called a \textit{subchain} (if we need to be specific, we say an \textit{$(a,b)$-subchain} or an \textit{$(a,c)$-subchain}). If $P_0, P_2, \dots, P_t$ are the $(a,b)$ subchains, and $Q_1, Q_3, \dots, Q_{t-1}$ are the $(a,c)$ subchains, appearing in that order on a meta-Kempe chain which has $u$ and $v$ as endpoints, then we denote the meta-Kempe chain as 
	$$u - P_0 - h_1 - Q_1 - h_2 - P_2 - h_3 - Q_3 - \dots - Q_{t-1} - h_t - P_t - v$$
	with $h_1, \dots, h_t$ being the switching vertices. 
	
	For a given minimal colouring $f$ of $G$, the edges coloured $0$ are called \emph{conflicting edges}. Each conflicting edge must have adjacent to it all three colours $1,2,3$ (otherwise the zero could be recoloured $1,2$ or $3$, contradicting the resistance). The set of conflicting edges can be partitioned into three sets: $A_1$, $A_2$ and $A_3$. $A_i$ is defined as the set of conflicting edges which have two of its four incident edges coloured $i$. The Parity Lemma is a well-known useful lemma which appears in many different forms in the literature on snarks. One such form, or consequence, is Lemma \ref{lem_parity}. The statements in Lemma \ref{lem_inducedcycles} have also appeared many times in the literature and are useful when dealing with snarks. These lemmas will be used implicitly in the proof of the main result.
	
	\begin{lemma} \label{lem_parity} {\rm \cite{Steffen1998}
			Let $G$ be a cubic graph. Let $f$ be a minimal colouring of $G$. Let $A = A_1 \cup A_2 \cup A_3$ be the set of conflicting edges. Then
			$$|A| \equiv |A_1| \equiv |A_2| \equiv |A_3| \mod 2.$$}  
	\end{lemma}
	
\begin{figure} [h!]	
		\begin{center}	
			\begin{tikzpicture}[
				x=1cm,
				y=1cm,
				font=\small,
				line width=0.8pt,
				line cap=round,
				line join=round,
				vertex/.style={
					circle, draw, fill=white,
					minimum size=5pt, inner sep=0pt
				},
				chain/.style={dash pattern=on 4pt off 4pt},
				colour/.style={fill=white, inner sep=2pt},
				explanation/.style={align=center, text width=4.8cm}
				]
				
				\foreach \shift/\leftcolour/\repeatcolour/\pair in {
					0/2/1/{2,3},
					6/1/2/{1,3}
				} {
					\begin{scope}[xshift=\shift cm]
						\coordinate (u) at (-0.65, 0);
						\coordinate (w) at ( 0.65, 0);
						
						\draw (u) -- node[colour, below=3pt] {$0$} (w);
						
						\draw (u) -- node[colour, above right] {$\leftcolour$}
						(-1.35, 0.70);
						\draw (w) -- node[colour, above left] {$3$}
						( 1.35, 0.70);
						
						\draw[chain]
						(-1.35, 0.70)
						.. controls (-2.15, 1.50) and (-1.90, 2.80) .. (0, 2.80)
						.. controls ( 1.90, 2.80) and ( 2.15, 1.50) .. (1.35, 0.70);
						\node at (0, 3.12) {$(\pair)$-Kempe chain};
						
						\draw (u) -- node[colour, below right] {$\repeatcolour$}
						(-1.35, -0.70);
						\draw (w) -- node[colour, below left] {$\repeatcolour$}
						( 1.35, -0.70);
						
						\node[vertex] at (u) {};
						\node[vertex] at (w) {};
					\end{scope}
				}
				
				\draw[decorate, decoration={brace, amplitude=7pt}]
				(-1.95, 3.65) -- (7.95, 3.65)
				node[midway, above=12pt, align=center]
				{These two induced cycles do not intersect.};
				
				\begin{scope}
					\coordinate (s) at (-0.65, -4.00);
					\coordinate (t) at ( 0.65, -4.00);
					
					\draw[chain]
					(-1.35, -0.70)
					.. controls (-2.00, -1.35) and (-2.00, -2.65) .. (-1.35, -3.30);
					\draw[chain]
					( 1.35, -0.70)
					.. controls ( 2.00, -1.35) and ( 2.00, -2.65) .. ( 1.35, -3.30);
					\node at (-1.25, -2.00) {$(1,3)$};
					\node at ( 1.25, -2.00) {$(1,2)$};
					
					\draw (-1.35, -3.30) -- node[colour, above right] {$3$} (s);
					\draw ( 1.35, -3.30) -- node[colour, above left]  {$2$} (t);
					\draw (s) -- node[colour, above=3pt] {$1$}
					node[below=4pt] {$e$} (t);
					
					\draw (s) -- node[colour, below right] {$2$} (-1.35, -4.70);
					\draw (t) -- node[colour, below left]  {$3$} ( 1.35, -4.70);
					\draw[chain] (-1.35, -4.70) -- (-1.95, -5.30);
					\draw[chain] ( 1.35, -4.70) -- ( 1.95, -5.30);
					
					\node[vertex] at (s) {};
					\node[vertex] at (t) {};
					
					\node[explanation, anchor=north] at (0, -6.25)
					{The $(1,3)$-chain and the $(1,2)$-chain\\
						meet at an edge $e$.};
				\end{scope}
				
				\begin{scope}[xshift=6cm]
					\coordinate (v) at (0, -4.00);
					\coordinate (z) at (0, -4.90);
					
					\draw[chain]
					(-1.35, -0.70)
					.. controls (-2.00, -1.35) and (-2.00, -2.65) .. (-0.75, -3.30);
					\draw[chain]
					( 1.35, -0.70)
					.. controls ( 2.00, -1.35) and ( 2.00, -2.65) .. ( 0.75, -3.30);
					\node at (-1.25, -2.00) {$(2,3)$};
					\node at ( 1.25, -2.00) {$(2,1)$};
					
					\draw (-0.75, -3.30) -- node[colour, below left] {$3$} (v);
					\draw ( 0.75, -3.30) -- node[colour, below right]  {$1$} (v);
					\draw (v) -- node[colour, right=3pt] {$2$} (z);
					
					\draw (z) -- node[colour, above left]  {$1$} (-0.70, -5.60);
					\draw (z) -- node[colour, above right] {$3$} ( 0.70, -5.60);
					\draw[chain] (-0.70, -5.60) -- (-1.15, -6.05);
					\draw[chain] ( 0.70, -5.60) -- ( 1.15, -6.05);
					
					\node[vertex, label=above:$v$] at (v) {};
					\node[vertex] at (z) {};
					
					\node[explanation, anchor=north] at (0, -6.25)
					{The $(2,3)$-chain and the $(2,1)$-chain\\
						meet at a vertex $v$.};
				\end{scope}
				
			\end{tikzpicture}
		\end{center}
		\caption{A diagrammatic representation of Lemma \ref{lem_inducedcycles}. In the diagram, we have as an example conflicting edges from $A_1$ and $A_2$ on the left and right, respectively. The induced cycles do not intersect. The $(a,b)$-chain and $(a,c)$-chain initiating at the conflicting edge intersect, and meet either at a vertex or at an edge.} \label{fig_lemma2}
	\end{figure}
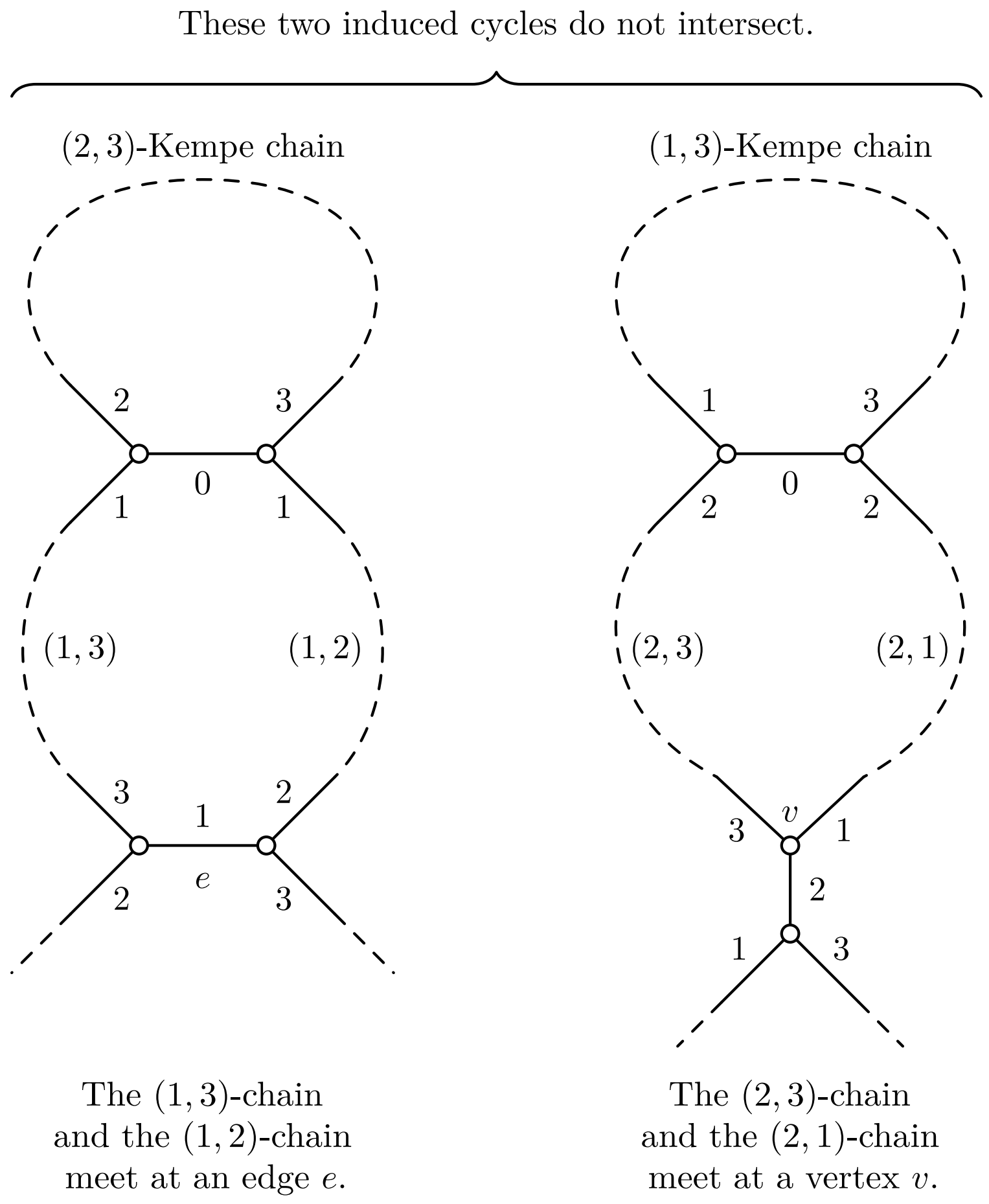

	\begin{lemma} \label{lem_inducedcycles} {\rm \cite{Steffen1998}
			Let $G$ be a cubic graph. Let $f$ be a minimal colouring of $G$. Let $e=uv \in A_a$ and let $u$ be incident to edges coloured $a,b$ and let $v$ be incident to edges coloured $a,c$, where $a,b,c$ are distinct colours in $\{1,2,3\}$. Then the following statements hold:
			\begin{enumerate}
				\item[(i)] The $(b,c)$-Kempe chain initiating at $u$ terminates at $v$. We say that this Kempe chain together with $e$ is a cycle \textit{induced} by $e$.
				\item[(ii)] The cycles induced by any two distinct conflicting edges do not intersect. 
				\item[(iii)] The $(a,c)$-Kempe chain initiating at $u$ and the $(a,b)$-Kempe chain initiating at $v$ must intersect. These chains either meet at a vertex or at an edge, as in Figure \ref{fig_lemma2}.  
		\end{enumerate}}   
	\end{lemma}

	\section{Main result}
	
	
	\begin{lemma} \label{lem_main} {\rm 
			Let $Z=\{u_1,\ldots,u_k,x_1,\ldots,x_k\}$ be a set of $2k$ distinct vertices and let $P_i$ be a $u_i$--$x_i$ path for each $i$ such that no vertex in $Z$ lies internally in any $P_i$. Let the paths $P_i$ be edge disjoint and let there be at most two paths intersecting at any vertex. Let $H = \bigcup_{i=1}^{k}P_i$. Then there exists $F$, a subgraph of $H$, such that $F$ has the following properties:
			\begin{itemize}
				\item There exists a set of edge-disjoint paths in $F$ whose endpoints are exactly the vertices in $Z$. 
				\item There are at most $2k-2$ intersecting vertices in $F$ (a vertex, $v$, is \textit{intersecting} in $F$ if it is incident to edges from $P_i$ and $P_j$ in $F$, for some $i \neq j$. If $v$ is incident to two edges from each of $P_i$ and $P_j$, \textit{then $v$ counts as two intersecting vertices}).  
				\item Each vertex in $F$ has degree two or four except that the vertices in $Z$ have degree one.
		\end{itemize}}
	\end{lemma}
	
	\begin{proof}
		By definition of $H$, every vertex in $Z$ has degree one, every intersecting vertex has degree four, and every other vertex has degree two.
		
		We first algorithmically determine a subforest of $H$, $T$, with at most $2k-2$ intersections and which contains $Z$, as follows. Initially, declare all $P_i$ unused. Now, choose an unused path $P_1$ and set $T_1=P_1$, and declare $P_1$ used. Then consider some unused $P_i$. If $P_i$ does not intersect $T_1$, leave it unused and consider a different path. If $P_i$ does intersect with $T_1$, traverse along $P_i$ from $u_i$ to $x_i$ and let $a$ and $b$ be the first and last intersecting vertices with $T_1$. Add to $T_1$ the $u_i$--$a$ and $b$--$x_i$ subpaths of $P_i$, and declare $P_i$ used. If $a=b$, then the entire $P_i$ is added, $a=b$ has degree four in $T_1$ and counts as two intersecting vertices. Continue doing this for each unused path. If at any point there are no unused paths which intersect $T_1$, then consider some unused path as another tree, $T_2$, and repeat the process. Continue this until all paths are used. Let the resultant graph be $T$. Note that each tree $T_i$ in $T$ contains an even number of vertices of $Z$. 
		
		It is clear that $T$ is a forest, because the initial path is a tree and every subpath added is a leaf segment. It is also clear that by construction, $T$ contains $Z$. Since the initial path adds no intersecting vertices to $T$, and each subsequent usage of a path adds at most $2$ intersecting vertices, and there are $k$ paths, it follows that the number of intersecting vertices in $T$ is at most $2(k-1)$. 
		
		There may exist vertices in $T$ with degree three. We aim to remove edges from $T$ such that the only vertices which remain have degree two or four, except for the vertices in $Z$ which have degree one. In each tree $T_j$, retain an edge, $e$, if and only if $T_j-e$ has an odd number of vertices of $Z$ in each of the two resulting components. Let $F$ be the forest induced by the retained edges. We claim that $F$ only contains vertices with degree two or four, except that the vertices in $Z$ have degree one.
		
		For any vertex $u_i$ in $Z$, the edge incident to $u_i$ must be retained, since removing it leaves a component with an odd number of vertices from $Z$ (the component containing just $u_i$). Therefore, each vertex in $Z$ has degree one in $F$.
		
		Let $v \notin Z, v \in T_j$, and let \(J_1,\ldots,J_d\) be the components of $T_j-v$, where $d=d_{T_j}(v)$. Let $e_i$ be the edge joining $v$ to $J_i$. Deleting $e_i$ separates $T_j$ into $J_i$ and its complement. Since $T_j$ contains an even number of vertices of $Z$, $J_i$ and its complement have the same parity of vertices of $Z$. Thus, $e_i$ is retained if and only if \(|V(J_i) \cap Z|\) is odd. Since 
		$$|V(J_1)\cap Z| + \dots + |V(J_d)\cap Z| = |V(T_j)\cap Z|,$$
		which is even, it must be that an even number of terms on the left side are odd. Therefore, an even number of the edges incident to $v$ are retained. It follows that $v$ has degree two or four in $F$.
		
				
		
		Since $T$ has at most $2k-2$ intersecting vertices, $F \subset T$ has at most $2k-2$ intersecting vertices.
		
		It remains to identify edge-disjoint paths pairing the vertices of $Z$. Now, since $F$ is a forest with $2k$ leaves (the vertices in $Z$), there must exist a path, $Q$, which initiates at one leaf and terminates at another. After removing $Q$, $F$ is a forest with $2k-2$ leaves (the remaining vertices in $Z$), so there must exist another path, $Q'$, which initiates at one leaf and terminates at another. And so on. It follows that vertices of $Z$ can be paired up by edge-disjoint paths in $F$. This completes the proof.
	\end{proof}
	
	\begin{figure} [h!]
		\begin{center}
			\begin{tikzpicture}[
				scale=0.8,
				x=1cm,
				y=0.85cm,
				font=\small,
				line cap=round,
				line join=round,
				rounded corners=9pt,
				original/.style={draw=black, line width=0.45pt},
				retained/.style={draw=black, line width=1.8pt},
				vertex/.style={
					circle, draw=black, fill=white,
					line width=0.6pt, minimum size=6pt, inner sep=0pt
				},
				endpoint/.style={vertex, minimum size=6.8pt},
				intersection/.style={vertex, line width=1.8pt, fill=black}
				]
				
				\coordinate (x1) at (0,  3.4);
				\coordinate (y1) at (10, 3.4);
				\coordinate (x2) at (0,  0);
				\coordinate (y2) at (10, 0);
				\coordinate (x3) at (0, -4.1);
				\coordinate (y3) at (10,-4.1);
				\coordinate (x4) at (0, -7.0);
				\coordinate (y4) at (10,-7.0);
				
				
				\def\PoneStart{
					(x1) -- (2.5, 3.4) -- (2.5, 0)
				}
				\def\PoneMiddle{
					(2.5, 0) -- (2.5,-8.2) -- (4.5,-8.2)
					-- (4.5, 2.0) -- (6.5, 2.0)
					-- (6.5,-1.8) -- (8.5,-1.8) -- (8.5, 0)
				}
				\def\PoneEnd{
					(8.5, 0) -- (8.5, 3.4) -- (y1)
				}
				
				\def\PthreeStart{
					(x3) -- (1.5,-4.1) -- (1.5, 0)
				}
				\def\PthreeMiddle{
					(1.5, 0) -- (1.5, 0.9) -- (3.5, 0.9)
					-- (3.5,-1.0) -- (5.5,-1.0)
					-- (5.5, 0.9) -- (7.5, 0.9) -- (7.5, 0)
				}
				\def\PthreeDescent{
					(7.5, 0) -- (7.5,-7.0)
				}
				\def\PthreeTurn{
					(7.5,-7.0) -- (7.5,-9.0) -- (8.5,-9.0) -- (8.5,-7.0)
				}
				\def\PthreeEnd{
					(8.5,-7.0) -- (8.5,-4.1) -- (y3)
				}
				
				
				\draw[original] \PoneStart;
				\draw[original] \PoneMiddle;
				\draw[original] \PoneEnd;
				
				\draw[original] (x2) -- (y2);
				
				\draw[original] \PthreeStart;
				\draw[original] \PthreeMiddle;
				\draw[original] \PthreeDescent;
				\draw[original] \PthreeTurn;
				\draw[original] \PthreeEnd;
				
				\draw[original] (x4) -- (y4);
				
				%
				%
				
				\draw[retained] \PoneStart;
				\draw[retained] (2.5, 0) -- (7.5, 0);
				\draw[retained] \PthreeDescent;
				\draw[retained] (7.5,-7.0) -- (x4);
				
				\draw[retained] (x2) -- (1.5, 0);
				\draw[retained] \PthreeStart;
				
				\draw[retained] \PoneEnd;
				\draw[retained] (8.5, 0) -- (y2);
				
				\draw[retained] \PthreeEnd;
				\draw[retained] (8.5,-7.0) -- (y4);
				
				
				\foreach \xx in {1.5,2.5,3.5,4.5,5.5,6.5,7.5,8.5} {
					\node[vertex] at (\xx,0) {};
				}
				
				\foreach \xx/\yy in {2.5/0.9,4.5/-1.0,6.5/0.9,7.5/-1.8} {
					\node[vertex] at (\xx,\yy) {};
				}
				
				\foreach \xx in {2.5,4.5,7.5,8.5} {
					\node[vertex] at (\xx,-7.0) {};
				}
				
				\foreach \xx/\yy in {1.5/0,2.5/0,7.5/0,8.5/0,7.5/-7.0,8.5/-7.0} {
					\node[intersection] at (\xx,\yy) {};
				}
				
				\foreach \i in {1,2,3,4} {
					\node[endpoint, label=left:$u_{\i}$] at (x\i) {};
					\node[endpoint, label=right:$x_{\i}$] at (y\i) {};
				}
				
				\node[above=5pt] at (0.75, 3.4) {$P_1$};
				\node[above=5pt] at (0.75, 0)   {$P_2$};
				\node[above=5pt] at (0.75,-4.1) {$P_3$};
				\node[above=5pt] at (0.75,-7.0) {$P_4$};
				
			\end{tikzpicture}
		\end{center}
		\caption{An example of the process described in Lemma \ref{lem_main}. We have $k=4$. Initially, set $T_1=P_2$, then use $P_1$, $P_3$ and $P_4$, in that order. Retaining edges as described in the proof results in $F$, which is represented by the thicker subpaths. The intersecting vertices are filled in black. The number of intersecting vertices respects the bound: $6 \leq 2(4) - 2$.} \label {fig_mainlemma}
	\end{figure}
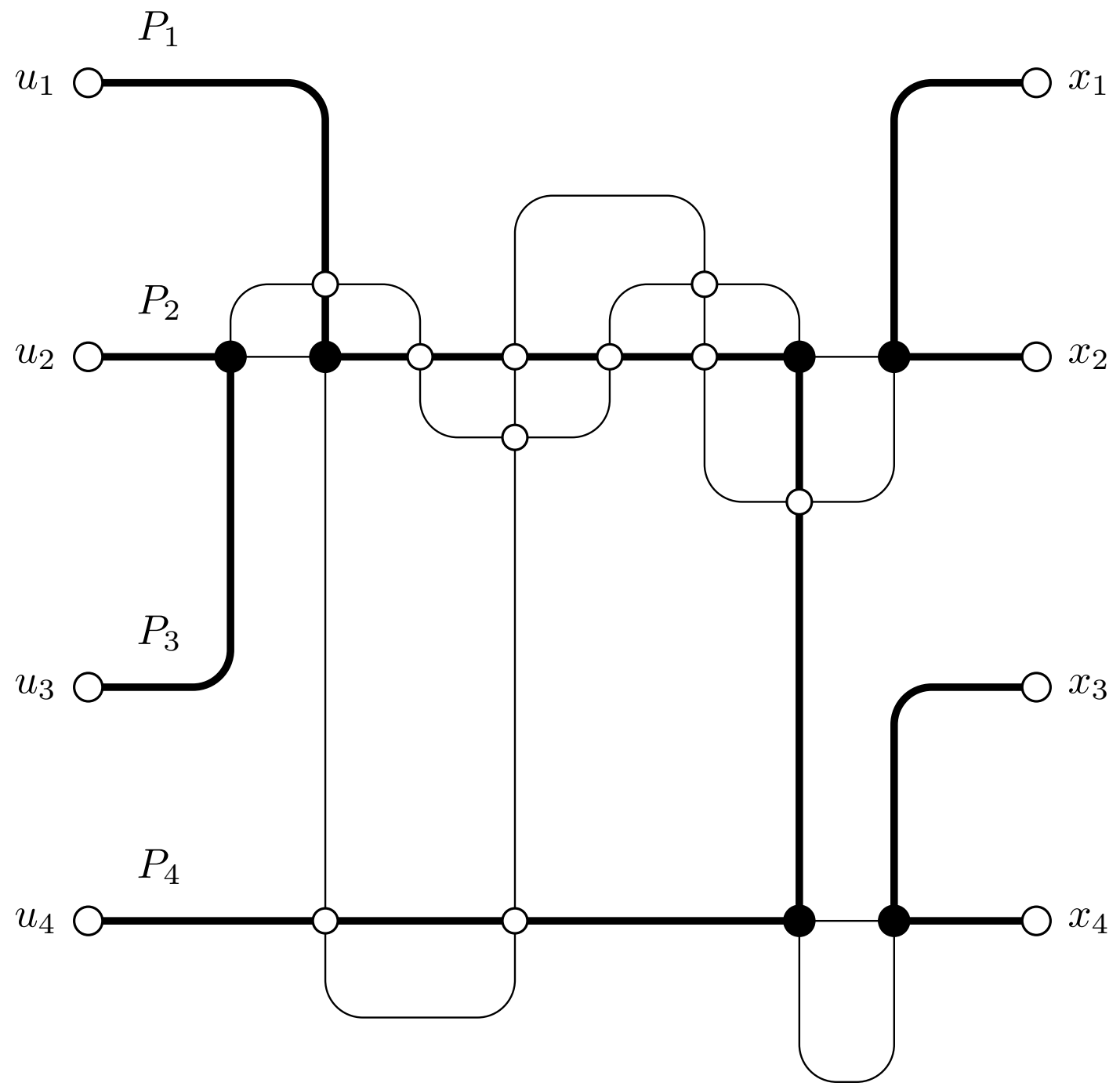
	
	
	\begin{theorem} \label{thm_main}
		Let $G$ be a snark. Then $\omega'(G) < 2r(G)$.
	\end{theorem}
	\begin{proof}
		Let $G$ be a bridgeless cubic graph with $r(G) = d$. Let $f$ be a minimal colouring of $G$. Let $|A_1| = p$, $|A_2| = q$ and $|A_3| = r$ so that $p+q+r=d$. Note that $p+q$ is even. We define the following notations: 
		\begin{itemize}
			\item Let $A_1 = \{u_1v_1, \dots, u_pv_p\}$ such that $u_i$ is incident to edges coloured $0,1,2$ and $v_i$ is incident to edges coloured $0,1,3$ for each $i$. 
			\item Let $A_2 = \{x_1y_1, \dots, x_qy_q\}$ such that $x_i$ is incident to edges coloured $0,2,1$ and $y_i$ is incident to edges coloured $0,2,3$ for each $i$.  
			\item Let $A_3 = \{w_1z_1, \dots, w_rz_r\}$ such that $w_i$ is incident to edges coloured $0,3,1$ and $z_i$ is incident to edges coloured $0,3,2$ for each $i$. 
			\item Let $Pu_i(1,3)$ be the maximal $(1,3)$-chain initiating at $u_i$.
			\item Let $Qx_i(2,3)$ be the maximal $(2,3)$-chain initiating at $x_i$.
			\item Let $Rw_i(3,2)$ be the maximal $(3,2)$-chain initiating at $w_i$. 
			\item Let $Rz_i(3,1)$ be the maximal $(3,1)$-chain initiating at $z_i$.
			\item Let $\mathcal{P} = \{Pu_i(1,3): 1 \leq i \leq p\}$. Note that any path $Pu_i(1,3) \in \mathcal{P}$ terminates at some $u_j$ or at some $z_j$ (in the former $Pu_i(1,3) = Pu_j(1,3)$, in the latter $Pu_i(1,3) = Rz_j(3,1)$).
			\item Let $\mathcal{Q} = \{Qx_i(2,3): 1 \leq i \leq q\}$. Note that any path $Qx_i(2,3) \in \mathcal{Q}$ terminates at some $x_j$ or at some $w_j$ (in the former $Qx_i(2,3) = Qx_j(2,3)$, in the latter $Qx_i(2,3) = Rw_j(3,2)$).
			\item Let $\mathcal{R} = \{Rw_i(3,2): 1 \leq i \leq r\} \cup \{Rz_i(3,1): 1 \leq i \leq r\}$. Note that any path $Rw_i(3,2) \in \mathcal{R}$ terminates at some $w_j$ or at some $x_j$ (in the former $Rw_i(3,2) = Rw_j(3,2)$, in the latter $Rw_i(3,2) = Qx_j(2,3)$), and $Rz_i(3,1) \in \mathcal{R}$ terminates at some $z_j$ or at some $u_j$ (in the former $Rz_i(3,1) = Rz_j(3,1)$, in the latter $Rz_i(3,1) = Pu_j(1,3)$). 
		\end{itemize}
		
		Note that $(1,3)$-chains in $\mathcal{P} \cup \mathcal{R}$ are either disjoint or equal, and $(2,3)$-chains in $\mathcal{Q} \cup \mathcal{R}$ are either disjoint or equal. 
	
		If $p=q=0$, then $r$ is even, and every vertex in $G$ has exactly one incident edge coloured $3$. The graph induced by the set of edges not coloured $3$ then induce a 2-factor. The odd components in the 2-factor are exactly the set of cycles induced by the conflicting edges. Therefore, $r(G) = \omega'(G) = r$, and the bound is satisfied. For the rest of the proof, we may assume that $p+q > 0$. 
		
		Let $Z = \{u_1, \dots, u_p, x_1, \dots, x_q\}$ and note its even cardinality. We aim to identify $(1,3)(2,3)$-meta-Kempe chains which pair up vertices in $Z$, each meta-Kempe chain having vertices from $Z$ as initial and terminal vertices, and containing no other vertices of $Z$. 
		
		We identify the $(1,3)(2,3)$-meta-Kempe chains as follows. If $Pu_{i_0}(1,3)$ equals some $Rz_{i_1}(3,1)$, then it must intersect with $Rw_{i_1}(3,2)$ at some vertex $h_1$ (such that $h_1$ is incident to an edge coloured $1$ from $Pu_{i_0}(1,3)$). Traverse along $Pu_{i_0}(1,3)$ from $u_{i_0}$ to $h_1$. Let $h_1$ be a switching vertex. If $Rw_{i_1}(3,2)$ equals some $Rw_{i_2}(3,2)$, then it must intersect with $Rz_{i_2}(3,1)$ at some vertex $h_2$ (such that $h_2$ is incident to an edge coloured $2$ from $Rw_{i_1}(3,2)$). Traverse from $h_1$ to $h_2$. Let $h_2$ be a switching vertex. If $Rz_{i_2}(3,1)$ equals some $Rz_{i_3}(3,1)$, then it must intersect with $Rw_{i_3}(3,2)$ at some vertex $h_3$ (such that $h_3$ is incident to an edge coloured $1$ from $Rz_{i_2}(3,1)$). Traverse from $h_2$ to $h_3$. Let $h_3$ be a switching vertex.
		Continue in this way, until some $Rz_{i_t}(3,1)$ equals some $Pu_{i_{t+1}}(1,3)$, or some $Rw_{i_t}(3,2)$ equals some $Qx_{i_{t+1}}(2,3)$. Then 
			$$u_{i_0} - Pu_{i_0}(1,3) - h_1-Rw_{i_1}(3,2) - h_2-Rz_{i_2}(3,1) - h_3 - \dots - h_t - Rz_{i_t}(3,1) - u_{i_{t+1}}$$
			$$\text{or}$$
			$$u_{i_0} - Pu_{i_0}(1,3) - h_1-Rw_{i_1}(3,2) - h_2-Rz_{i_2}(3,1) - h_3 - \dots - h_t - Rw_{i_t}(3,2) -  x_{i_{t+1}}$$
		is the meta-Kempe chain which pairs up $u_{i_0}$ and $u_{i_{t+1}}$, or $u_{i_0}$ and $x_{i_{t+1}}$, in $Z$.
		It may be that $t=0$. This simply means that $Pu_{i_0}(1,3)$ equals some $Pu_{i_{1}}(1,3)$. 
		
		If the initial Kempe chain is $Qx_{i_0}(1,3)$ instead of $Pu_{i_0}(1,3)$, we proceed in a similar way to get 
		$$x_{i_0} - Qx_{i_0}(2,3) - h_1 - Rz_{i_1}(3,1) - h_2 - Rw_{i_2}(3,2) - h_3 - \dots - h_t - Rz_{i_t}(1,3) - u_{i_{t+1}}$$
		$$\text{or}$$
		$$x_{i_0} - Qx_{i_0}(2,3) - h_1 - Rz_{i_1}(3,1) - h_2 - Rw_{i_2}(3,2) - h_3 - \dots - h_t - Rw_{i_t}(2,3) - x_{i_{t+1}}$$
		as the meta-Kempe chain which pairs up $x_{i_0}$ and $u_{i_{t+1}}$, or $x_{i_0}$ and $x_{i_{t+1}}$, in $Z$. Again, it may be that $t=0$. This simply means that the initial $Qx_{i_0}(2,3)$ equals some $Qx_{i_{1}}(2,3)$. 
		
		Since the $(1,3)$ chains in $\mathcal{P} \cup \mathcal{R}$ are disjoint or equal, and the $(2,3)$ chains in $\mathcal{Q} \cup \mathcal{R}$ are disjoint or equal, it follows that there are no duplicate vertices in $\{u_{i_0}, w_{i_1}, z_{i_1}, \dots, w_{i_t}, z_{i_t}, u_{i_{t+1}}\} \cup \{u_{j_0}, w_{j_1}, z_{j_1}, \dots, w_{j_s}, z_{j_s}, u_{j_{s+1}}\}$, or similarly for any pair of meta-Kempe chains. Therefore, the meta-Kempe chain pairs up vertices in $Z$, and the set of pairs is a partition of $Z$. See Figure \ref{fig_metakempechains}.  
		
\begin{figure} [h!]
	\begin{center}	
		\begin{tikzpicture}[
			line width=0.6pt,
			vertex/.style={circle,draw,fill=white,inner sep=0pt,minimum size=5pt},
			chain/.style={dashed,line width=0.6pt},
			bold chain/.style={dashed,line width=2pt},
			colour/.style={font=\small,inner sep=1pt},
			name/.style={font=\scriptsize,inner sep=1pt},
			every node/.style={font=\small}
			]
			
			\begin{scope}[yshift=4.5cm]
				\foreach \x/\k/\a/\b/\c/\d in {2/i/2/1/3/1,9/j/1/2/1/3} {
					\draw (\x,0) -- node[colour,left=2pt] {$0$} (\x,-0.9);
					\draw (\x,0) -- ++(-0.5,0.5) node[colour,pos=0.25,left=4pt] {$\a$};
					\draw (\x,0) -- ++( 0.5,0.5) node[colour,pos=0.25,right=4pt] {$\b$};
					\draw (\x,-0.9) -- ++(-0.5,-0.5) node[colour,pos=0.25,left=4pt] {$\c$};
					\draw (\x,-0.9) -- ++( 0.5,-0.5) node[colour,pos=0.25,right=4pt] {$\d$};
					\node[vertex] at (\x,0) {};
					\node[vertex] at (\x,-0.9) {};
					\node[name,above=6pt] at (\x,0) {$u_\k$};
					\node[name,below=6pt] at (\x,-0.9) {$v_\k$};
				}
				\draw[bold chain] (2.5,0.5)
				.. controls (3.2,1.2) and (7.8,1.2) .. (8.5,0.5);
				\node[above=5pt] at (5.5,0) {$Pu_i(1,3) = Pu_j(1,3)$};
			\end{scope}
			
			\begin{scope}
				\foreach \x/\k/\a/\b/\c/\d in {
					0/0/2/1/3/1, 11/{t+1}/1/2/1/3
				} {
					\draw (\x,1.3) -- node[colour,left=2pt] {$0$} (\x,0.4);
					\draw (\x,1.3) -- ++(-0.5,0.5) node[colour,pos=0.25,left=4pt] {$\a$};
					\draw (\x,1.3) -- ++( 0.5,0.5) node[colour,pos=0.25,right=4pt] {$\b$};
					\draw (\x,0.4) -- ++(-0.5,-0.5) node[colour,pos=0.25,left=4pt] {$\c$};
					\draw (\x,0.4) -- ++( 0.5,-0.5) node[colour,pos=0.25,right=4pt] {$\d$};
					\node[vertex] at (\x,1.3) {};
					\node[vertex] at (\x,0.4) {};
					\node[name,above=12pt] at (\x,1.3) {$u_{i_{\k}}$};
					\node[name,below=12pt] at (\x,0.4) {$v_{i_{\k}}$};
				}
				
				\draw[bold chain] (0.5,1.8) .. controls (1.5,2) and (0.5,0) ..  (1.4,0);
				\draw[bold chain] (2.2,0) -- (4.2,0);
				\draw[bold chain] (5,0) -- (6.3,0);
				\node at (7.1,0) {\textbf{$\cdots$}};
				\draw[bold chain] (7.9,0) -- (8.8,0);
				\draw[bold chain] (9.6,0)
				.. controls (10.5,0) and (9.5,2) .. (10.5,1.8);
				%
				
				\node at (2,1.3) {$Pu_{i_0}(1,3)$};
				\node at (9,1.3) {$Rz_{i_{t}}(3,1)$};
				\node[above=2pt] at (3.2,0) {$Rw_{i_1}(3,2)$};
				\node[above=2pt] at (5.8,0) {$Rz_{i_2}(3,1)$};
				\node[above=2pt] at (8.1,0) {$Rw_{i_{t-1}}(3,2)$};
				
				\foreach \x/\k/\L/\R/\a/\b in {
					1.8/1/z_{i_1}/w_{i_1}/2/1,
					4.6/2/w_{i_2}/z_{i_2}/1/2,
					9.2/t/w_{i_t}/z_{i_t}/1/2
				} {
					
					\draw[chain] (\x-0.4,-0.3)
					.. controls (\x-1.25,-0.85) and (\x-1.25,-1.25)
					.. (\x-1.05,-1.55);
					\draw[chain] (\x+0.4,-0.3)
					.. controls (\x+1.25,-0.85) and (\x+1.25,-1.25)
					.. (\x+1.05,-1.55);
					
					\draw (\x-0.55,-2.05) -- (\x-1.05,-1.55)
					node[colour,pos=0.65,right=3pt] {$3$};
					\draw (\x+0.55,-2.05) -- (\x+1.05,-1.55)
					node[colour,pos=0.65,left=3pt] {$3$};
					\draw (\x-0.55,-2.05) -- (\x+0.55,-2.05)
					node[colour,midway,below=2pt] {$0$};
					
					\node[vertex, scale=3] at (\x,0) {};
					\node[name,below=10pt] at (\x,0) {$h_\k$};
				}
				
				\draw (1.8-0.55,-2.05) -- ++(-0.5,-0.5)
				node[colour,pos=0.65,right=3pt] {$1$};
				\draw (1.8+0.55,-2.05) -- ++(0.5,-0.5)
				node[colour,pos=0.65,left=3pt] {$2$};
				\draw (4.6-0.55,-2.05) -- ++(-0.5,-0.5)
				node[colour,pos=0.65,right=3pt] {$2$};
				\draw (4.6+0.55,-2.05) -- ++(0.5,-0.5)
				node[colour,pos=0.65,left=3pt] {$1$};
				\draw (9.2-0.55,-2.05) -- ++(-0.5,-0.5)
				node[colour,pos=0.65,right=3pt] {$2$};
				\draw (9.2+0.55,-2.05) -- ++(0.5,-0.5)
				node[colour,pos=0.65,left=3pt] {$1$};
				
				\node[vertex] at (1.8-0.55,-2.05) {};
				\node[vertex] at (1.8+0.55,-2.05) {};
				\node[vertex] at (4.6-0.55,-2.05) {};
				\node[vertex] at (4.6+0.55,-2.05) {};
				\node[vertex] at (9.2-0.55,-2.05) {};
				\node[vertex] at (9.2+0.55,-2.05) {};
				
				\node[name,left=8pt] at (1.8-0.55,-2.05) {$w_{i_1}$};
				\node[name,right=8pt] at (1.8+0.55,-2.05) {$z_{i_1}$};
				\node[name,left=8pt] at (4.6-0.55,-2.05) {$z_{i_2}$};
				\node[name,right=8pt] at (4.6+0.55,-2.05) {$w_{i_2}$};
				\node[name,left=8pt] at (9.2-0.55,-2.05) {$z_{i_t}$};
				\node[name,right=8pt] at (9.2+0.55,-2.05) {$w_{i_t}$};
				
				\node at (7.1,-2.05) {$\cdots$};
			\end{scope}

			\begin{scope} [yshift=-5.5cm]
				\foreach \x/\k/\a/\b/\c/\d in {
					0/0/2/1/3/1, 11/{t+1}/2/1/2/3
				} {
					\draw (\x,1.3) -- node[colour,left=2pt] {$0$} (\x,0.4);
					\draw (\x,1.3) -- ++(-0.5,0.5) node[colour,pos=0.25,left=4pt] {$\a$};
					\draw (\x,1.3) -- ++( 0.5,0.5) node[colour,pos=0.25,right=4pt] {$\b$};
					\draw (\x,0.4) -- ++(-0.5,-0.5) node[colour,pos=0.25,left=4pt] {$\c$};
					\draw (\x,0.4) -- ++( 0.5,-0.5) node[colour,pos=0.25,right=4pt] {$\d$};
					\node[vertex] at (\x,1.3) {};
					\node[vertex] at (\x,0.4) {};
				}
				
				\node[name,above=12pt] at (0,1.3) {$u_{i_0}$};
				\node[name,above=12pt] at (11,1.3) {$x_{i_{t+1}}$};
				\node[name,below=12pt] at (0,0.4) {$v_{i_0}$};
				\node[name,below=12pt] at (11,0.4) {$y_{i_{t+1}}$};
				
				\draw[bold chain] (0.5,1.8) .. controls (1.5,2) and (0.5,0) ..  (1.4,0);
				\draw[bold chain] (2.2,0) -- (4.2,0);
				\draw[bold chain] (5,0) -- (6.3,0);
				\node at (7.1,0) {\textbf{$\cdots$}};
				\draw[bold chain] (7.9,0) -- (8.8,0);
				\draw[bold chain] (9.6,0)
				.. controls (10.5,0) and (9.5,2) .. (10.5,1.8);

				\node at (2,1.3) {$Pu_{i_0}(1,3)$};
				\node at (9,1.3) {$Rw_{i_{t}}(3,2)$};
				\node[above=2pt] at (3.2,0) {$Rw_{i_1}(3,2)$};
				\node[above=2pt] at (5.8,0) {$Rz_{i_2}(3,1)$};
				\node[above=2pt] at (8.1,0) {$Rz_{i_{t-1}}(3,1)$};
				
				\foreach \x/\k/\L/\R/\a/\b in {
					1.8/1/z_{i_1}/w_{i_1}/2/1,
					4.6/2/w_{i_2}/z_{i_2}/1/2,
					9.2/t/w_{i_t}/z_{i_t}/1/2
				} {
					
					\draw[chain] (\x-0.4,-0.3)
					.. controls (\x-1.25,-0.85) and (\x-1.25,-1.25)
					.. (\x-1.05,-1.55);
					\draw[chain] (\x+0.4,-0.3)
					.. controls (\x+1.25,-0.85) and (\x+1.25,-1.25)
					.. (\x+1.05,-1.55);
					
					\draw (\x-0.55,-2.05) -- (\x-1.05,-1.55)
					node[colour,pos=0.65,right=3pt] {$3$};
					\draw (\x+0.55,-2.05) -- (\x+1.05,-1.55)
					node[colour,pos=0.65,left=3pt] {$3$};
					\draw (\x-0.55,-2.05) -- (\x+0.55,-2.05)
					node[colour,midway,below=2pt] {$0$};
					
					\node[vertex, scale=3] at (\x,0) {};
					\node[name,below=10pt] at (\x,0) {$h_\k$};
				}
				
				\draw (1.8-0.55,-2.05) -- ++(-0.5,-0.5)
				node[colour,pos=0.65,right=3pt] {$1$};
				\draw (1.8+0.55,-2.05) -- ++(0.5,-0.5)
				node[colour,pos=0.65,left=3pt] {$2$};
				\draw (4.6-0.55,-2.05) -- ++(-0.5,-0.5)
				node[colour,pos=0.65,right=3pt] {$2$};
				\draw (4.6+0.55,-2.05) -- ++(0.5,-0.5)
				node[colour,pos=0.65,left=3pt] {$1$};
				\draw (9.2-0.55,-2.05) -- ++(-0.5,-0.5)
				node[colour,pos=0.65,right=3pt] {$1$};
				\draw (9.2+0.55,-2.05) -- ++(0.5,-0.5)
				node[colour,pos=0.65,left=3pt] {$2$};
				
				\node[vertex] at (1.8-0.55,-2.05) {};
				\node[vertex] at (1.8+0.55,-2.05) {};
				\node[vertex] at (4.6-0.55,-2.05) {};
				\node[vertex] at (4.6+0.55,-2.05) {};
				\node[vertex] at (9.2-0.55,-2.05) {};
				\node[vertex] at (9.2+0.55,-2.05) {};

				\node[name,left=8pt] at (1.8-0.55,-2.05) {$w_{i_1}$};
				\node[name,right=8pt] at (1.8+0.55,-2.05) {$z_{i_1}$};
				\node[name,left=8pt] at (4.6-0.55,-2.05) {$z_{i_2}$};
				\node[name,right=8pt] at (4.6+0.55,-2.05) {$w_{i_2}$};
				\node[name,left=8pt] at (9.2-0.55,-2.05) {$w_{i_t}$};
				\node[name,right=8pt] at (9.2+0.55,-2.05) {$z_{i_t}$};

				\node at (7.1,-2.05) {$\cdots$};
			\end{scope}
		\end{tikzpicture}
	\end{center}
	\caption{Diagrammatic representations of meta-Kempe chains in snark $G$ with a minimal colouring $f$. The thicker dashed lines represent the meta-Kempe chains. Circles labeled $h_i$ represent switching vertices. } \label{fig_metakempechains}		
\end{figure}
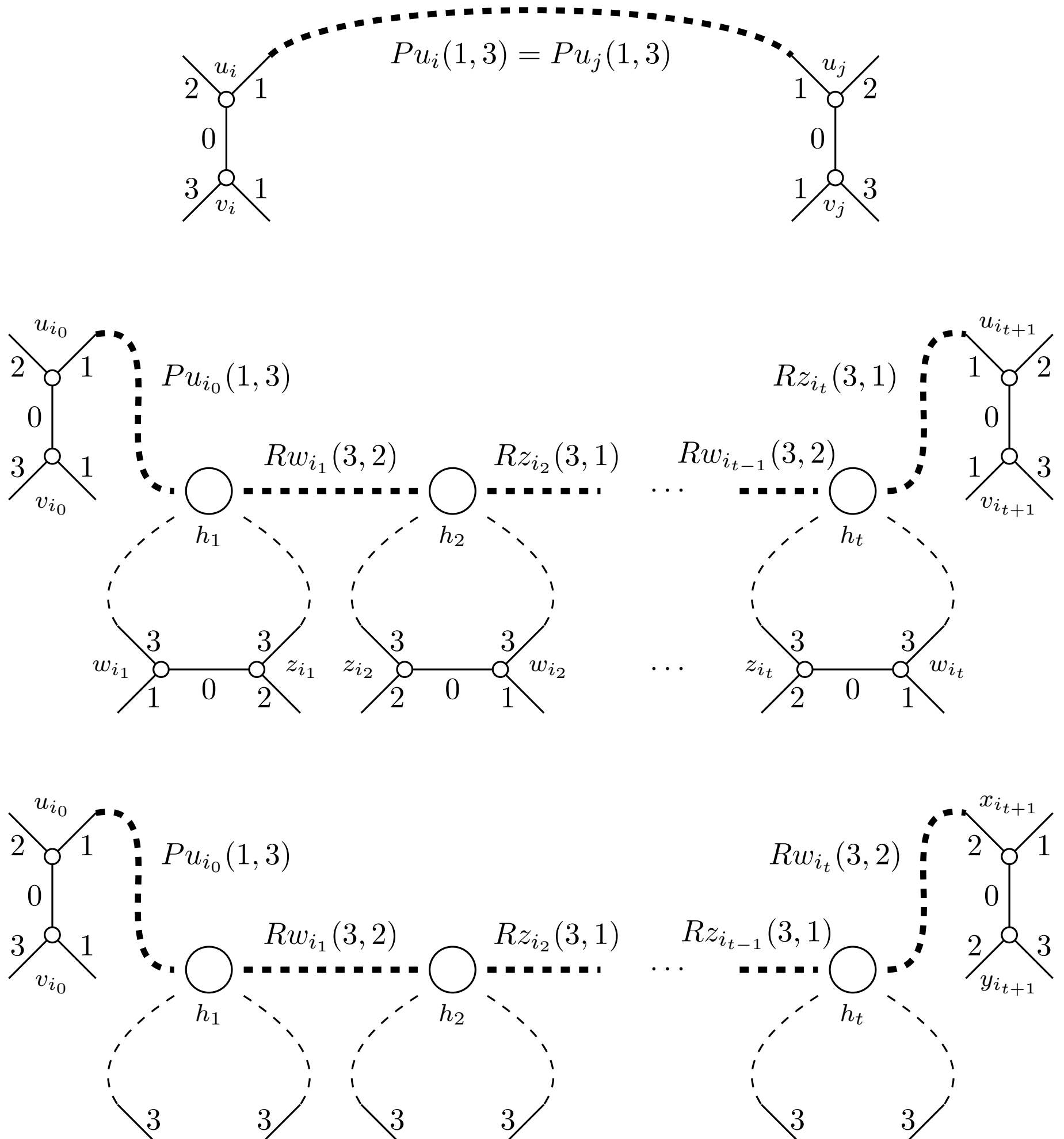		
		
		It may be that the identified $(1,3)(2,3)$-meta-Kempe chain intersects with itself. In such a case, we find a subpath which is also a meta-Kempe chain as follows. Let 
		$$u_0 - Pu_{i_0}(1,3) - h_1-Rw_{i_1}(3,2) - h_2-Rz_{i_2}(3,1) - h_3 - \dots - h_t - Rz_{i_t}(3,1) - u_{i_{t+1}}$$
		be a meta-Kempe chain which intersects with itself. Firstly, we note that any intersection must be a $(1,3)$ subchain intersecting with a $(2,3)$ subchain. Consider the first subchain $Pu_{i_0}(1,3)$. If it does not intersect with any other subchain (except at $h_1$), then we proceed to consider the next subchain, $Rw_{i_1}(2,3)$. If $Pu_{i_0}(1,3)$ does intersect with some other $Rw_{i_s}(2,3)$ subchain before encountering $h_1$, then consider the first such intersecting vertex, say $h_s$ (such that $h_s$ is incident to an edge coloured $1$ from $Pu_{i_0}(1,3)$). Let $h_s$ be a switching vertex instead. Traverse along $Rw_{i_s}(2,3)$ from $h_s$ toward $u_{i_{t+1}}$ to form a new meta-Kempe chain. Note that the new meta-Kempe chain has gained a switching vertex ($h_s$), but also has lost at least the switching vertex $h_1$. In the new chain, the initial subchain $Pu_{i_0}(1,3)$ now does not intersect with any further $(2,3)$ subchain. Now, in the new meta-Kempe chain, we consider the second subchain (which is a $(2,3)$ subchain) and conduct the same process. Continue in this way until no subchain intersects with any other subchain (except at the switching vertices), and the meta-Kempe chain is a path. Note, the number of switching vertices in the resultant chain is less than or equal to the number in the original meta-Kempe chain.  
		
		
		Let $p+q=2k$. Relabel vertices of $Z$ to $Z=\{u_1, \dots, u_k, x_1, \dots, x_k\}$ and relabel the meta-Kempe chains to $P_1, \dots, P_k$ such that $P_i$ is the meta-Kempe chain with $u_i$ and $x_i$ as endpoints. Let $H = \bigcup_{i=1}^{k}P_i$. Contract all $3$-coloured intersecting edges between all pairs $P_i$ and $P_j$ in $H$. $H$ now satisfies the conditions of Lemma \ref{lem_main}. Therefore, there exists $F \subset H$ which has the properties listed in Lemma \ref{lem_main}. Now reverse the contraction of each edge which was contracted in $H$. If the reverse-contracted edge is adjacent to two or four edges in $F$, then retain it in $F$. Otherwise, ignore it.
		
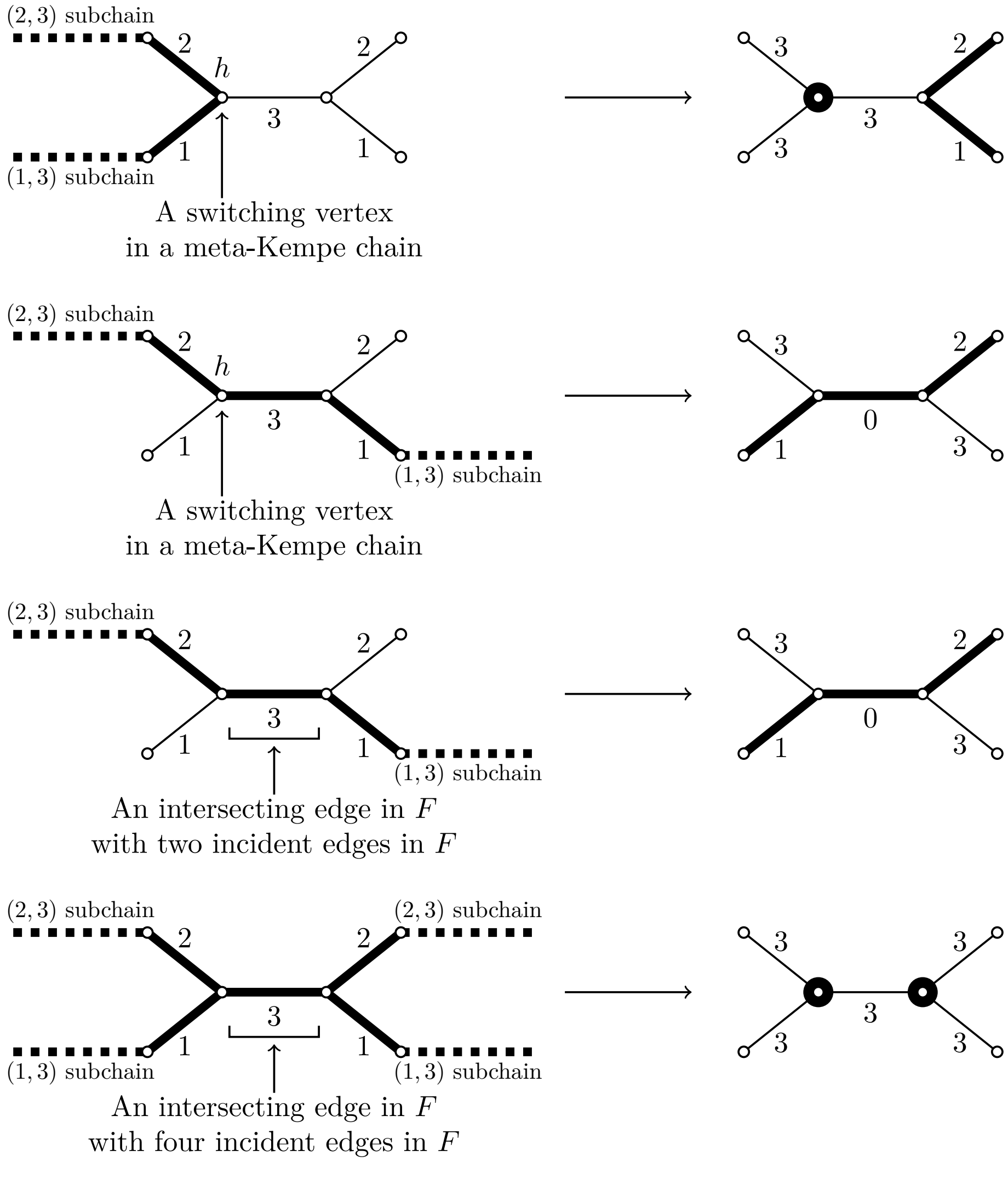
\begin{figure} [h!]
	\begin{center}
		\begin{tikzpicture}[
			scale=0.9,
			x=1cm, y=1cm,
			line width=0.7pt,
			every node/.style={font=\small}
			]
			
			
			\begin{scope}
				\draw[line width=3pt] (-1,0.8) -- node[above] {$2$} (0,0);
				\draw[line width=3pt] (-1,-0.8) -- node[below] {$1$} (0,0);
				\node[name,above=3pt] at (0,0) {$h$};
				\draw (0,0) -- node[below] {$3$} (1.4,0);
				\draw (1.4,0) -- node[above] {$2$} (2.4,0.8);
				\draw (1.4,0) -- node[below] {$1$} (2.4,-0.8);
				
				\draw[dashed, line width=3pt] (-2.8,0.8) -- (-1,0.8);
				\draw[dashed, line width=3pt] (-2.8,-0.8) -- (-1,-0.8);
				\node[above,scale=0.8] at (-1.9,0.8) {$(2,3)$ subchain};
				\node[below,scale=0.8] at (-1.9,-0.8) {$(1,3)$ subchain};
				
				\foreach \p in {(-1,0.8),(-1,-0.8),(0,0),
					(1.4,0),(2.4,0.8),(2.4,-0.8)}
				\draw[fill=white] \p circle (2pt);
				
				\node[align=center] at (0.7,-1.8)
				{A switching vertex \\ in a meta-Kempe chain};
				\draw[->] (0,-1.35) -- (0,-0.2);
			\end{scope}
			
			\begin{scope}[xshift=8cm]
				\draw (-1,0.8) -- node[above] {$3$} (0,0);
				\draw (-1,-0.8) -- node[below] {$3$} (0,0) node[circle, fill] {};
				\draw (0,0) -- node[below] {$3$} (1.4,0);
				\draw[line width=3pt] (1.4,0) -- node[above] {$2$} (2.4,0.8);
				\draw[line width=3pt] (1.4,0) -- node[below] {$1$} (2.4,-0.8);
				
				\foreach \p in {(-1,0.8),(-1,-0.8),(0,0),
					(1.4,0),(2.4,0.8),(2.4,-0.8)}
				\draw[fill=white] \p circle (2pt);
			\end{scope}
			
			\begin{scope}[yshift=-4cm]
				\draw[line width=3pt] (-1,0.8) -- node[above] {$2$} (0,0);
				\draw (-1,-0.8) -- node[below] {$1$} (0,0);
				\node[name,above=3pt] at (0,0) {$h$};
				\draw[line width=3pt] (0,0) -- node[below] {$3$} (1.4,0);
				\draw (1.4,0) -- node[above] {$2$} (2.4,0.8);
				\draw[line width=3pt] (1.4,0) -- node[below] {$1$} (2.4,-0.8);
				
				\draw[dashed, line width=3pt] (-2.8,0.8) -- (-1,0.8);
				\draw[dashed, line width=3pt] (2.4,-0.8) -- (4.2,-0.8);
				\node[above,scale=0.8] at (-1.9,0.8) {$(2,3)$ subchain};
				\node[below,scale=0.8] at (3.3,-0.8) {$(1,3)$ subchain};
				
				\foreach \p in {(-1,0.8),(-1,-0.8),(0,0),
					(1.4,0),(2.4,0.8),(2.4,-0.8)}
				\draw[fill=white] \p circle (2pt);

				\node[align=center] at (0.7,-1.8)
				{A switching vertex \\ in a meta-Kempe chain};
				\draw[->] (0,-1.35) -- (0,-0.2);
				
			\end{scope}
			
			\begin{scope}[xshift=8cm,yshift=-4cm]
				\draw (-1,0.8) -- node[above] {$3$} (0,0);
				\draw[line width=3pt] (-1,-0.8) -- node[below] {$1$} (0,0);
				\draw[line width=3pt] (0,0) -- node[below] {$0$} (1.4,0);
				\draw[line width=3pt] (1.4,0) -- node[above] {$2$} (2.4,0.8);
				\draw (1.4,0) -- node[below] {$3$} (2.4,-0.8);
				
				\foreach \p in {(-1,0.8),(-1,-0.8),(0,0),
					(1.4,0),(2.4,0.8),(2.4,-0.8)}
				\draw[fill=white] \p circle (2pt);
			\end{scope}

			\begin{scope}[yshift=-8cm]
				\draw[line width=3pt] (-1,0.8) -- node[above] {$2$} (0,0);
				\draw (-1,-0.8) -- node[below] {$1$} (0,0);
				\draw[line width=3pt] (0,0) -- node[below] {$3$} (1.4,0);
				\draw (1.4,0) -- node[above] {$2$} (2.4,0.8);
				\draw[line width=3pt] (1.4,0) -- node[below] {$1$} (2.4,-0.8);
				
				\draw[dashed, line width=3pt] (-2.8,0.8) -- (-1,0.8);
				\draw[dashed, line width=3pt] (2.4,-0.8) -- (4.2,-0.8);
				
				\node[above,scale=0.8] at (-1.9,0.8) {$(2,3)$ subchain};
				\node[below,scale=0.8] at (3.3,-0.8) {$(1,3)$ subchain};
				
				\foreach \p in {(-1,0.8),(-1,-0.8),(0,0),
					(1.4,0),(2.4,0.8),(2.4,-0.8)}
				\draw[fill=white] \p circle (2pt);
				
				\node[align=center] at (0.7,-1.8)
				{An intersecting edge in $F$ \\ with two incident edges in $F$};
				\draw (0.1,-0.45) -- (0.1,-0.6) -- (1.3,-0.6) -- (1.3,-0.45);
				\draw[->] (0.7,-1.35) -- (0.7,-0.7);
				
			\end{scope}
			
			\begin{scope}[xshift=8cm,yshift=-8cm]
				\draw (-1,0.8) -- node[above] {$3$} (0,0);
				\draw[line width=3pt] (-1,-0.8) -- node[below] {$1$} (0,0);
				\draw[line width=3pt] (0,0) -- node[below] {$0$} (1.4,0);
				\draw[line width=3pt] (1.4,0) -- node[above] {$2$} (2.4,0.8);
				\draw (1.4,0) -- node[below] {$3$} (2.4,-0.8);
				
				\foreach \p in {(-1,0.8),(-1,-0.8),(0,0),
					(1.4,0),(2.4,0.8),(2.4,-0.8)}
				\draw[fill=white] \p circle (2pt);
			\end{scope}
			
			\begin{scope}[yshift=-12cm]
				\draw[line width=3pt] (-1,0.8) -- node[above] {$2$} (0,0);
				\draw[line width=3pt] (-1,-0.8) -- node[below] {$1$} (0,0);
				\draw[line width=3pt] (0,0) -- node[below] {$3$} (1.4,0);
				\draw[line width=3pt] (1.4,0) -- node[above] {$2$} (2.4,0.8);
				\draw[line width=3pt] (1.4,0) -- node[below] {$1$} (2.4,-0.8);
				
				\draw[dashed,line width=3pt] (-2.8,0.8) -- (-1,0.8);
				\draw[dashed,line width=3pt] (-2.8,-0.8) -- (-1,-0.8);
				\draw[dashed,line width=3pt] (2.4,0.8) -- (4.2,0.8);
				\draw[dashed,line width=3pt] (2.4,-0.8) -- (4.2,-0.8);
				
				\node[above,scale=0.8] at (-1.9,0.8) {$(2,3)$ subchain};
				\node[below,scale=0.8] at (-1.9,-0.8) {$(1,3)$ subchain};
				\node[above,scale=0.8] at (3.3,0.8) {$(2,3)$ subchain};
				\node[below,scale=0.8] at (3.3,-0.8) {$(1,3)$ subchain};
				
				\foreach \p in {(-1,0.8),(-1,-0.8),(0,0),
					(1.4,0),(2.4,0.8),(2.4,-0.8)}
				\draw[fill=white] \p circle (2pt);
				
				\node[align=center] at (0.7,-1.8)
				{An intersecting edge in $F$ \\ with four incident edges in $F$};
				\draw (0.1,-0.45) -- (0.1,-0.6) -- (1.3,-0.6) -- (1.3,-0.45);
				\draw[->] (0.7,-1.35) -- (0.7,-0.7);
			\end{scope}
			
			\begin{scope}[xshift=8cm,yshift=-12cm]
				\draw (-1,0.8) -- node[above] {$3$} (0,0);
				\draw (-1,-0.8) -- node[below] {$3$} (0,0) node[circle, fill] {};
				\draw (0,0) -- node[below] {$3$} (1.4,0) node[circle, fill] {};
				\draw (1.4,0) -- node[above] {$3$} (2.4,0.8);
				\draw (1.4,0) -- node[below] {$3$} (2.4,-0.8);
				
				\foreach \p in {(-1,0.8),(-1,-0.8),(0,0),
					(1.4,0),(2.4,0.8),(2.4,-0.8)}
				\draw[fill=white] \p circle (2pt);
			\end{scope}

			\draw[->] (4.6,0) -- (6.3,0);
			\draw[->] (4.6,-4) -- (6.3,-4);
			\draw[->] (4.6,-8) -- (6.3,-8);
			\draw[->] (4.6,-12) -- (6.3,-12);
			
		\end{tikzpicture}
	\end{center}
	\caption{Diagrammatic representation of colour swapping on the meta-Kempe chains and the resulting even factor. The top two diagrams show both examples, without loss of generality, of how colour swapping occurs at switching vertices in the meta-Kempe chains. The bottom two diagrams show both examples of how two meta-Kempe chains intersect in $F$. The thicker edges on the left represent the meta-Kempe chains before colour swapping. The thicker edges and vertices on the right represent the even factor after colour swapping and subsequent removal of all $3$-coloured edges.} \label{fig_swappingcolours}
\end{figure}

		On each $(1,3)$ and $(2,3)$ subchain in $F$, we now swap colours $1$ and $3$, and $2$ and $3$, respectively, except for $3$-coloured edges in $F$ which are adjacent to at least one $1$-coloured edge in $F$ at one end and at least one $2$-coloured edge in $F$ at the other end. Some of these $3$-coloured edges will be left coloured $3$, and some will be coloured $0$ instead, as follows. 
		\begin{itemize}
			\item At switching vertices in $F$ (from the meta-Kempe chains):
			\begin{itemize}
				\item[$\triangleright$] if the vertex was not incident to a $3$-coloured edge in $F$, it now has all three incident edges coloured $3$ in $G$.
				\item[$\triangleright$] if the vertex was incident to a $3$-coloured edge in $F$, colour the edge $0$ instead, and notice that it now has two adjacent edges coloured $3$, one at either end.
			\end{itemize} 
			\item At reverse-contracted intersecting edges in $F$ (which must be $3$-coloured):
			\begin{itemize}
				\item[$\triangleright$] If it has two adjacent edges in $F$, colour it $0$ instead, and notice that it now has two adjacent edges coloured $3$, one at either end.
				\item[$\triangleright$] If it has four adjacent edges in $F$, leave the intersecting edge coloured $3$, and notice that it now has all four adjacent edges also coloured $3$. Recall from Lemma \ref{lem_main} that when such an edge was contracted, the identified vertex counted as two intersecting vertices in $F$.
			\end{itemize}
		\end{itemize}
		Every vertex in $G$ is now incident with exactly one edge or exactly three edges coloured $3$. Therefore, removing the $3$-coloured edges from $G$ leaves an even factor of $G$. See Figure \ref{fig_swappingcolours}. The only odd components in this even factor are the isolated vertices, and possibly the cycles which contain a $0$-coloured edge. Therefore, we have that each switching vertex from the meta-Kempe chains, each intersecting edge in $F$, and each initial conflicting edge with regard to $f$, adds at most one odd component to the even factor. 
		
		Note that the resistance of $G$ is $r + 2k$. For the $k$ meta-Kempe chains $P_1, \dots, P_k$, let $t_1, \dots, t_k$ be the number of switching vertices on each, respectively. Also, in $F$ there are at most $(2k-2)$ intersecting edges. Note as well that there exists an injective mapping from the set of switching vertices to the set of edges in $A_3$, implying that $(t_1 + \dots + t_k) \leq r$. Therefore, we have that 
		\begin{align*}
			\frac{\omega'(G)}{r(G)} 
			& \leq \frac{(t_1 + \dots + t_k) + (2k-2) + (r + 2k)}{(r + 2k)} \\
			& \leq \frac{(r) + (2k-2) + (r + 2k)}{(r + 2k)} \\
			& = \frac{2(r + 2k) - 2 }{(r + 2k)} \\
			& < 2.
		\end{align*}
	\end{proof}	
		

	As it turns out, oddness and weak oddness are remarkably different.  

	\begin{problem} {\rm
		Is the bound $\omega'(G) < 2r(G)$ asymptotically sharp?}
	\end{problem}
	
	\section*{Acknowledgment}
	
	We acknowledge the use of AI ChatGPT Plus in the production of this manuscript. Our aim here is to best describe all instances where AI assisted meaningfully in the mathematics. Other instances where AI assisted include only the writing of the bibliography and production of diagrams in Figures (not in the conceptualization of diagrams, only in the typesetting. Pictures were drawn by hand and AI was asked to produce LaTeX code). All novel mathematics outside of Lemma 3 is original work conceived and produced entirely by the author (although engaging with the AI in a conversational manner about the problem may have assisted the author in attaining clarity of ideas, much like discussion with a peer would). That is to say, the problem was simplified to the formulation of Lemma 3 without any direct AI assistance. Once the problem was reduced in this way, we decided to write up the formulation of Lemma 3 and prompt AI to provide a proof, after brief investigation which suggested that the statement was true. It took a few iterations to get the formulation precisely correct. The AI did not contribute in refining the formulation. It was just made clear to the author by the AI's response that further precision was required. Except in one case where it noted an implicit assumption which was then added. That being, ``... such that no vertex in $Z$ lies internally in any $P_i$". The proof of Lemma 3 provided by the AI was dense, contained some `human'-like errors such as making an implicit assumption, and was not rigorous as per the standards of the discipline, in our opinion. Ultimately, by completely understanding and refining the proof ourselves to our satisfaction (as we would typically with our own produced proofs), we take ownership of the proof of Lemma 3 while acknowledging AI's assistance, as we tentatively put forth our current belief that this is the best practice in this new age of AI-assisted mathematics.

\end{document}